\documentclass[12pt,a4paper,draft]{article}

\usepackage[english]{babel}
\usepackage{amssymb}
\usepackage{latexsym,amsfonts,amsmath,longtable,mathrsfs}
\usepackage[unicode,final,hyperindex]{hyperref}
\usepackage{amsthm}

\DeclareMathOperator{\GL}{GL} \DeclareMathOperator{\Aut}{Aut}
\DeclareMathOperator{\SL}{SL}
\DeclareMathOperator{\Sym}{Sym}

\renewcommand{\le}{\leqslant}

\begin{document}

\newtheorem{lem}{Lemma}
\newtheorem{thm}[lem]{Theorem}
\newtheorem{cor}[lem]{Corollary}
\newtheorem{prop}[lem]{Proposition}
\theoremstyle{remark}
\newtheorem{exmp}[lem]{{\slshape Example}}

{\bfseries 2000MSC} 20D10

\begin{center}
{\bfseries An existence criterion for Hall subgroups of finite
groups\footnote{The work is supported by  RFBR, projects 08-01-00322, 10-01-00391, and 10-01-90007,  ADTP ``Development of the Scientific
Potential of
Higher
School'' of the Russian Federal Agency for Education (Grant
2.1.1.419), and Federal Target Grant "Scientific and
educational personnel of innovation Russia" for 2009-2013 (government
contract No. 02.740.11.0429). The second author gratefully acknowledges the support from
Deligne 2004 Balzan prize in mathematics. }}

Danila O. Revin\footnote{Institute of Mathematics, pr-t Acad. Koptyug, 4, Novosibirsk, Russia,
630090; e-mail: revin@math.nsc.ru}

Evgeny P. Vdovin\footnote{Corresponding author; Institute of Mathematics, pr-t Acad. Koptyug, 4, Novosibirsk, Russia,
630090; e-mail: vdovin@math.nsc.ru}
\end{center}

\begin{small}
{\bfseries Abstract.} We obtain an existence criterion for Hall subgroups of finite
groups in terms of a composition series. As a corollary we provide a solution to Problem 5.65 from the Kourovka notebook.

\end{small}

\section{Introduction}
\label{Introduction}

Let  $\pi$ be a set of primes. We denote by  $\pi'$ the set of all primes not in
$\pi$, by $\pi(n)$ the set of prime divisors of a positive integer $n$,
for a finite group $G$ we denote $\pi(|G|)$ by $\pi(G)$. A positive integer
$n$ with $\pi(n)\subseteq\pi$ is called a $\pi$-number,  a group $G$ with
$\pi(G)\subseteq \pi$ is called a $\pi$-group. A subgroup $H$ of $G$ is called
a {\it $\pi$-Hall subgroup}, if $\pi(H)\subseteq\pi$ and $\pi(|G:H|)\subseteq
\pi'$.  According to \cite{Hall} we  say that $G$ {\it satisfies
$E_\pi$} (or briefly $G\in E_\pi$), if $G$ contains a $\pi$-Hall subgroup. If
$G\in E_\pi$ and every two $\pi$-Hall subgroups are conjugate, then we say
that  $G$ {\it satisfies $C_\pi$} ($G\in C_\pi$). If $G\in C_\pi$ and
each $\pi$-subgroup of $G$ is included in a  $\pi$-Hall
subgroup of $G$, then we  say that $G$ {\it satisfies $D_\pi$} ($G\in
D_\pi$). A group satisfying $E_\pi$ (resp. $C_\pi$, $D_\pi$) is also called
an $E_\pi$-group
(resp. a $C_\pi$-group, a $D_\pi$-group).
Let $A,B,H$ be subgroups of $G$ such that $B\unlhd A$ and $H\leq
G$. Then $N_H(A/B)=N_H(A)\cap N_H(B)$ is the {\em normalizer} of $A/B$ in $H$. If $x\in
N_H(A/B)$, then $x$ induces an automorphism of $A/B$ by $Ba\mapsto B x^{-1}ax$.
Thus there exists a homomorphism $N_H(A/B)\rightarrow \Aut(A/B)$. The image of
$N_H(A/B)$ under this homomorphism is denoted by $\Aut_H(A/B)$ and is called a
{\em group of induced automorphisms} of $A/B$, while the kernel of this
homomorphism is denoted by~$C_H(A/B)$ and is called the {\em centralizer} of $A/B$ in $H$.

All results of the paper depend on the classification of finite
simple groups.  We can avoid its use with the following definition. A finite group is said to be a {\em $K$-group}, if all of its
composition factors are known simple groups, i.e., are either cyclic of prime order, or alternating, or of Lie type, or sporadic. If we
replace the term ``finite group'' by the term ``finite $K$-group'', then all results remain valid without the use of the
classification of finite simple groups.

In \cite[Theorem~7.7]{RevVdoDpiFinal} we proved the following theorem.

\begin{thm}\label{DpiCriterion}
Let $A$ be a normal subgroup of a finite group $G$. Then $G\in D_\pi$ if and
only if $A\in D_\pi$ and~${G/A\in D_\pi}$.
\end{thm}

By using this theorem we obtain that a finite group $G$ satisfies $D_\pi$ if
and only if every  composition factor satisfies $D_\pi$. In \cite{RevDpiDAN} an arithmetic description of finite simple
$D_\pi$-gro\-ups was obtained; thus there exists a precise exhaustive arithmetic criterion for
determining whether a finite group $G$ satisfies~$D_\pi$. Our aim is to find a similar criterion for~$E_\pi$.

The following easy proposition is well known  (see \cite[Lemma~1]{Hall}).

\begin{prop} \label{base}
Let $A$ be a normal subgroup of $G$.  If $H$ is a $\pi$-Hall subgroup of $G$, then
$H \cap A$ is a $\pi$-Hall subgroup of $A$, and $HA/A$ is a
$\pi$-Hall subgroup of~$G/A$.
\end{prop}

It follows from this that a normal subgroup and a factor group of an $E_\pi$-gro\-up satisfy $E_\pi$. But  an extension of an
$E_\pi$-gro\-up by an $E_\pi$-gro\-up may fail to satisfy $E_\pi$ (see Example \ref{ExtensionIsNotEpi}). Thus a
criterion for  a finite group to satisfy $E_\pi$ will be different in nature to a  criterion
for~$D_\pi$.

In \cite[Theorem~3.5]{GrossExistence} Gross proved
the following theorem.

\begin{thm}\label{GrossExistenceCriterion}
Let $1=G_0<G_1<\ldots< G_n=G$ be a composition series of a finite group $G$
which is a refinement of a chief series of $G$. Then the following are
equivalent:
\begin{itemize}
\item[{\em (a)}] $H\in E_\pi$ for all subgroups $H$ such that
$H^{(\infty)}$ is subnormal in $G$, where $H^{(\infty)}$ is the intersection of all members of the derived series of~$H$.
\item[{\em (b)}] $\Aut_G(A/A^\ast)\in E_\pi$ for all $A\in
\mathcal{A}(G)$.
\item[{\em (c)}] $\Aut_G(G_i/G_{i-1})\in E_\pi$ for all $i$ with $1\le i\le n$.
\item[{\em (d)}] $\Aut_G(H/K)\in E_\pi$ for all composition factors $H/K$
of~$G$.
\end{itemize}
\end{thm}

The
definition of $\Aut_G(A/B)$ is given above; the symbol $\mathcal{A}(G)$ denotes the set of all atoms of $G$. Recall that a subgroup $A$ of
$G$ is called an {\em atom}, if $A$ is subnormal in $G$, $A=A'$, and $A$ has exactly one maximal normal subgroup. 
If in statement (a) we take $H=G$, then we obtain that $G\in E_\pi$. Thus the
condition $\Aut_G(G_i/G_{i-1})\in E_\pi$ for all $i$ with $1\le i\le n$ implies
that $G\in E_\pi$.
The following result is our main theorem.

\begin{thm}\label{ExistenceCriterion} 
Let $1=G_0<G_1<\ldots< G_n=G$ be a composition series of a finite group $G$.
If, for some $i$,
$\Aut_G(G_i/G_{i-1})\not\in E_\pi$, then~${G\not\in E_\pi}$.
\end{thm}

The recently completed classification of $\pi$-Hall subgroups in finite simple groups (see \cite{RevVdoArxive}) is the main technical
tool for the proof of the theorem.
Combining Theorems \ref{GrossExistenceCriterion} and \ref{ExistenceCriterion}
we obtain a criterion for a finite group to satisfy~$E_\pi$.

\begin{cor}\label{GrossExistenceCriterionImproved}
Let $1=G_0<G_1<\ldots< G_n=G$ be a composition series of a finite group $G$
which is a refinement of a chief series of $G$. Then the following are
equivalent:
\begin{itemize}
\item[{\em (a)}] $H\in E_\pi$ for all subgroups $H$ such that
$H^{(\infty)}$ is subnormal in $G$, where $H^{(\infty)}$ is the intersection of all members of the derived series of~$H$.
\item[{\em (b)}] $\Aut_G(A/A^\ast)\in E_\pi$ for all $A\in
\mathcal{A}(G)$.
\item[{\em (c)}] $\Aut_G(G_i/G_{i-1})\in E_\pi$ for all $i$ with $1\le i\le n$.
\item[{\em (d)}] $\Aut_G(H/K)\in E_\pi$ for all composition factors $H/K$
of~$G$.
\item[{\em (e)}] $G\in E_\pi$.
\end{itemize}
\end{cor}

If $G_i/G_{i-1}$ is cyclic, then $\Aut_G(G_i/G_{i-1})\leq \Aut(G_i/G_{i-1})$ is
also cyclic, hence $\Aut_G(G_i/G_{i-1})\in E_\pi$. Thus we need to check
that $\Aut_G(G_i/G_{i-1})\in E_\pi$  only for non-abelian composition factors.
Moreover, let $$1=G_0<G_1<G_2<\ldots<G_k=G$$ be a refinement of a chief series
$$1=G_0=G_{i_0}< G_{i_1}<\ldots< G_{i_m}=G_k=G.$$ Then
$G_{i_{j+1}}/G_{i_{j}}=T_1\times\ldots\times T_s$ is a direct product of
isomorphic sim\-p\-le groups. Assume that $G_{i_{j+1}}/G_{i_{j}}$ is non-abelian;
thus $T_1\simeq T_2\simeq\ldots\simeq T_s$ are non-abelian finite simple groups.
We may choose the numbering of $T_1,\ldots,T_s$ so that $G_{i_{j}+1}$ is a
complete preimage of $T_1$ in $G$, $G_{i_{j}+2}$ is a complete preimage of
$T_1\times T_2$ in $G$, etc. Since for every $l\le s$ we have
$G_{i_{j}}\leq C_G(G_{i_{j}+l}/G_{i_{j}+l-1})$, we have by
\cite[Lemma~1.2]{Vdo} that $\Aut_G(G_{i_{j}+l}/G_{i_{j}+l-1})\simeq
\Aut_{G/G_{i_{j}}}((G_{i_{j}+l}/G_{i_j})/(G_{i_{j}+l-1}/G_{i_j}))$.
Assume that $j=0$ so
$T_1,\ldots,T_s$ are subgroups of $G$. Since $T_1\times\ldots\times T_s$ is a
minimal normal subgroup of $G$, then $G$ acts transitively by conjugation on
the set $\{T_1,\ldots,T_s\}$ and $N_G(T_1),\ldots,\linebreak N_G(T_s)$ are all conjugate
in $G$. By \cite[Theorem~3.3.10]{Rob} we obtain
that $$N_G((T_1\times\ldots\times T_l)/(T_1\times\ldots\times T_{l-1}))\leq
N_G(T_l)$$ and 
\begin{multline*}
C_G((T_1\times\ldots\times T_l)/(T_1\times\ldots\times T_{l-1}))\\ =C_G(T_l)\cap N_G((T_1\times\ldots\times
T_l)/(T_1\times\ldots\times T_{l-1})) 
\end{multline*}
for each $l\le s$; hence  $\Aut_G(G_{i_{0}+l}/G_{i_{0}+l-1})\leq
\Aut_G(G_{i_{0}+1}/G_{i_{0}})$. Similarly, in the general case
$\Aut_G(G_{i_{j}+l}/G_{i_{j}+l-1})\leq
\Aut_G(G_{i_{j}+1}/G_{i_{j}})$. Now \cite[Corollary~3.3]{GrossExistence}
implies that if $\Aut_G(G_{i_{j}+1}/G_{i_{j}})\in E_\pi$, then $\Aut_G(G_{i_{j}+l}/G_{i_{j}+l-1})\in E_\pi$ for each $l\le
s$. Thus we obtain
the following

\begin{cor}\label{InducedAutomorphismsEpi}
Let $1=G_0<G_1<G_2<\ldots<G_k=G$ be a refinement of a chief series
$1=G_0=G_{i_0}< G_{i_1}<\ldots< G_{i_m}=G_k=G$. Then  $G\in E_\pi$ if and only if for every non-abelian
$G_{i_{j+1}}/G_{i_j}$ we have  $\Aut_G(G_{i_j+1}/G_{i_j})\in E_\pi$.
\end{cor}

So we have the following problem: describe the almost simple $E_\pi$-groups. Recall that a finite group $G$ is called almost simple if
generalized Fitting subgroup $F^\ast(G)$ of $G$ is simple, i.e., if $S\leq G\leq \Aut(S)$ for a non-abelian finite simple group $S$. If an
arithmetic description of finite almost simple
$E_\pi$-gro\-ups is obtained, then Corollary \ref{InducedAutomorphismsEpi} would give
an exhaustive arithmetic criterion for determining whether a finite group $G$ satisfies~$E_\pi$.

At the end of the paper we prove some corollaries to the main result. Corollary \ref{SubDirectProduct} is the main result of \cite{Ved}, but
the solution presented in  \cite{Ved} is known to contain a gap\footnote{In the proof of Lemma 3 the equality $H\cap
\mathrm{Aut}(P_i)=\mathrm{Aut}_H(P_i)$ in line 16 of p.220 is incorrect, so the statement $\Aut_H(P_i)\in E_\pi$ is not proven.}, which at
the time of writing has yet to be filled.

\begin{cor}\label{SubDirectProduct}
A finite subdirect product of $E_\pi$-gro\-ups satisfies $E_\pi$.
\end{cor}

Combining Proposition \ref{base} and Corollary \ref{SubDirectProduct} one can obtain an affirmative answer to~\cite[Problem~5.65]{kour} and
\cite[Problem~18]{Shemetkov}.

\begin{cor}\label{FormationOfEpiGroups}
For every set of primes $\pi$ the  class of all $E_\pi$-gro\-ups is a formation.
\end{cor}

In view of Proposition \ref{base} a homomorphic image of a $\pi$-Hall subgroup is a $\pi$-Hall subgroup (in the homomorphic image of the
$E_\pi$-gro\-up). The following corollary shows that all $\pi$-Hall subgroups of a homomorphic image can be obtained in this way.

\begin{cor}\label{HomomorphicImage}
Every $\pi$-Hall of a homomorphic image of an $E_\pi$-gro\-up $G$ is the image of a $\pi$-Hall subgroup of~$G$.
\end{cor}

In contrast, there are examples showing that in general a normal subgroup of an $E_\pi$-gro\-up may possess $\pi$-Hall subgroups that do
not lie in
any $\pi$-Hall subgroup of the whole group (cf. Lemma~\ref{HallExist} and Example~\ref{ex2}).

The main technical tool in the proof of Theorem \ref{ExistenceCriterion} is the
following lemma, which is of independent interest.

\begin{lem}\label{EpiCyclic}
Let $S$ be a non-abelian finite simple $E_\pi$-gro\-up and suppose that
$S< G\leq
\Aut(S)$ and $G\not\in E_\pi$. Then $2,3\in\pi\cap \pi(S)$ and there exists an element
$x$ of $G$ such that $\langle x,S\rangle\not\in E_\pi$ and the order of $x$ is either a power of $2$, or a power of~$3$.
\end{lem}

To prove this lemma we shall use the theorem about the number of classes of conjugate $\pi$-Hall subgroups in finite simple groups~\cite[Theorem~1.1]{RevVdoArxive}.

\section{Notation and preliminary results}

The term ``group'' always means a ``finite group''.
By $\pi$ we always denote a set of primes. The expressions $H\leq G$ and $H\unlhd G$ mean respectively
that $H$ is a subgroup and a normal subgroup of $G$. The symbol $K_\pi(G)$ denotes the set of
classes of conjugate $\pi$-Hall subgroups of~$G$ and $k_\pi(G)=\vert K_\pi(G)\vert$. The following lemma follows from the Schur-Zassenhaus
Theorem (see
\cite[Theorems~D6 and~D7]{Hall} or \cite[Chapter~5, Theorem~3.7]{Suz}, for example).

\begin{lem}\label{NormalSeriesWithpiprimefactors}
If every factor of a subnormal series of $G$ are either a $\pi$-group or a
$\pi'$-group, then~${G\in D_{\pi}}$.
\end{lem}

\begin{lem}\label{InducedAutomorphism}
Let $M$, $N$ be normal subgroups of a group $G$ such that $M\cap N=1$, and $A,B$ subgroups of $M$ such that $B\unlhd A$.  Let
$$\overline{\phantom{G}}:G\rightarrow G/N=\overline{G}$$ be the natural homomorphism. Then $N\leq
C_G(A/B)$, $\overline{N_G(A/B)}=N_{\overline{G}}(\overline{A}/\overline{B})$, $\overline{C_G(A/B)}=
C_{\overline{G}}(\overline{A}/\overline{B})$, and
$\Aut_G(A/B)\simeq \Aut_{\overline{G}}(\overline{A}/\overline{B})$.
\end{lem}

{\bfseries Proof.}\ \ The inclusion $N\leq C_G(A/B)$ is true since the product of  $M$ and $N$ is direct. The inclusions
$\overline{N_G(A/B)}\subseteq N_{\overline{G}}(\overline{A}/\overline{B})$ and $\overline{C_G(A/B)}\subseteq
C_{\overline{G}}(\overline{A}/\overline{B})$ are evident. We need to show the reverse inclusions in order to complete
the proof.

Let $\overline{x}\in N_{\overline{G}}(\overline{A}/\overline{B})$ for some $x\in G$. Then $A^xN=AN$ and $B^xN=BN$. Since
$M\trianglelefteq G$, then $A^x\leq M$ and $B^x\leq M$. Applying the coordinate projection map $M\times N\rightarrow M$ to both parts of the
equalities $A^xN=AN$ and $B^xN=BN$, we obtain
$A^x=A$ and $B^x=B$. It follows that $x\in N_G(A/B)$, and consequently we obtain the inclusion  $\overline{N_G(A/B)}\supseteq
N_{\overline{G}}(\overline{A}/\overline{B})$.

Let $\overline{x}\in C_{\overline{G}}(\overline{A}/\overline{B})$. Then for every $a\in A$ there exist $b\in B$, $g\in N$ such that
$a^x=bag$. Since $a^x\in A$ and $A,B\leq M$, the 
projection  $M\times N\rightarrow M$ with kernel $N$ leaves  $a^x$ stable, and at the same time maps it into $ba$. Therefore $a^x=ba$.  Thus
 $\overline{C_{G}(A/B)}\supseteq C_{\overline{G}}(\overline{A}/\overline{B})$.\qed

The following lemma is known and its proof can be found in \cite[Lem\-ma~2.1(e)]{RevVdoArxive}, for example.

\begin{lem}\label{HallExist}
Let $A$ be a normal $E_\pi$-sub\-gro\-up of $G$ such that $G/A$ is a
$\pi$-gro\-up, and $M$ a $\pi$-Hall subgroup of $A$. Then there exists  a $\pi$-Hall subgroup
$H$ of $G$ satisfying $H\cap A=M$  if and only if the set~${\{M^a\mid a\in A\}}$ is closed under conjugation by~$G$.
\end{lem}

An extension of an $E_\pi$-gro\-up by an $E_\pi$-gro\-up may fail to
possess a $\pi$-Hall subgroup, as is shown in the following known example.

\begin{exmp}\label{ExtensionIsNotEpi}
Let $\pi=\{2,3\}$, $G=\GL_3(2)=\SL_3(2)$ be
a
group of order $168=2^3\cdot 3\cdot 7$. From
\cite[Theorem  1.2]{RevHallp} 
it follows that $G$ has exactly
two classes of $\pi$-Hall subgroups with representatives
$$
\left(
\begin{array}{c@{}c}
\fbox{
$\begin{array}{c}
\!\!\!
\GL_2(2)
\!\!\! \\
\end{array}$
}
& *\\
0 &\fbox{1}
\end{array}
\right)\text{ and }
\left(
\begin{array}{c@{}c}
\fbox{1}& *\\
 0&\fbox{$\begin{array}{c}
\!
\GL_2(2)
\!
\\
\end{array}$}
\end{array}
\right).
$$
The first  consists of  line stabilizers in the natural
representation of $G$, and the second  consists of plain
stabilizers. The map $\iota : x\in G\mapsto (x^t)^{-1}$, where $x^t$
denotes the transpose of the matrix  $x$, is an automorphism of order $2$
of  $G$. It interchanges classes of $\pi$-Hall subgroups, hence by Proposition \ref{base} and Lemma \ref{HallExist} the natural  extension
$\widehat{G}=G\rtimes\langle\iota\rangle$ does not possess a $\pi$-Hall
subgroup.
\end{exmp}

The following example shows that  a normal subgroup of an $E_\pi$-gro\-up may possess
$\pi$-Hall subgroups that are not contained in
a $\pi$-Hall subgroup of the whole group.

\begin{exmp}\label{ex2}
Let $\pi=\{2,3\}$. Let $G=\GL_5(2)=\SL_5(2)$ be a
group of order $99999360=2^{10}\cdot 3^2\cdot 5\cdot 7 \cdot
31$. Let $\iota : x\in G\mapsto (x^t)^{-1}$ and $\widehat{G}=G\rtimes
\langle\iota\rangle$ be a natural semidirect product. From
\cite[Theorem~1.2]{RevHallp} it follows that there exist $\pi$-Hall subgroups
of $G$, and every such a subgroup is a stabilizer of a series of subspaces
$V=V_0<V_1<V_2<V_3=V$, where $V$ is the natural module of $G$, and $\dim
V_k/V_{k-1}\in \{1,2\}$ for  $k=1,2,3$. Therefore, there are
three conjugacy classes of $\pi$-Hall subgroups of $G$ with representatives
$$
H_1=\left(
\begin{array}{c@{}c@{}c}
\fbox{$\begin{array}{c}
\\
\!
\GL_2(2)
\!
\\
\\
\end{array}$}& &*\\
 &\fbox{1}& \\
0& &
\fbox{$\begin{array}{c}
\\
\!
\GL_2(2)
\!
\\
\\
\end{array}$}
\end{array}
\right),
$$
$$H_2=
{\left(
\begin{array}{c@{}c@{}c}
\fbox{1}& &*\\
 &\fbox{$\begin{array}{c}
\\
\!
\GL_2(2)
\!
\\
\\
\end{array}$}& \\
0& &
\fbox{$\begin{array}{c}
\\
\!
\GL_2(2)
\!
\\
\\
\end{array}$}
\end{array}
\right),
\text{ and }
H_3=\left(
\begin{array}{c@{}c@{}c}
\fbox{$\begin{array}{c}
\\
\!
\GL_2(2)
\!
\\
\\
\end{array}$}& &*\\
 &\fbox{$\begin{array}{c}
\\
\!
\GL_2(2)
\!
\\
\\
\end{array}$}& \\
0& &
\fbox{1}
\end{array}
\right).}
$$
Note that  $N_G(H_k)=H_k$, $k=1,2,3$, since  $H_k$ is parabolic. The class
containing $H_1$ is $\iota$-invariant. So the Frattini argument implies that
$\widehat{G}=GN_{\widehat{G}}(H_1)$, whence $|N_{\widehat{G}}(H_1):N_G(H_1)|=2$
and  $N_{\widehat{G}}(H_1)$ is a $\pi$-Hall subgroup of  $\widehat{G}$. Moreover 
$\iota$ interchanges classes containing  $H_2$ and $H_3$. So, as in the
previous example, these subgroups are not contained in $\pi$-Hall subgroups of~$\widehat{G}$.
\end{exmp}

\section{Proof of the Main Theorem and Corollaries}

\noindent {\bfseries Proof of Lemma \ref{EpiCyclic}.} Suppose that $S$ is a finite simple $E_\pi$-gro\-up and $G$ is  chosen so that $S<
G\leq
\Aut(S)$ and $G\not\in E_\pi$. Assume also that  for every proper subgroup $M$ of $G$ with $S\leq M$ we have $M\in E_\pi$, i.e., $G$ is a
minimal subgroup of $\Aut(S)$ containing $S$ subject to~${G\not\in E_\pi}$.

Notice that $G/S$ is a $\pi$-gro\-up. Indeed, by the Schreier conjecture
$G/S\leq \Aut(S)/S$ is solvable, hence it satisfies
$D_\pi$. Consider a $\pi$-Hall subgroup $M/S$ of $G/S$. If $G/S$ is not a $\pi$-gro\-up, then $M$ is a proper subgroup of $G$,  hence it
satisfies~$E_\pi$ in view of the minimality of~$G$. So there exists a $\pi$-Hall subgroup $H$ of~$M$. Since $\vert G:M\vert=\vert
(G/S):(M/S)\vert$ and $M/S$ is a $\pi$-Hall subgroup of $G/S$, we obtain that $\vert G:M\vert$ is a $\pi'$-num\-ber, so $H$ is a $\pi$-Hall
subgroup of~$G$. This contradicts~${G\not\in E_\pi}$.

We have that $G$, acting by conjugation, permutes  elements
of $K_\pi(S)$, so we obtain a homomorphism $\varphi:G\rightarrow \Sym(K_\pi(S))\simeq \Sym_{k_\pi(S)}$. By Lemma \ref{HallExist} we obtain
that $G\not\in E_\pi$ if and only if $G$ does not leave invariant any conjugacy class of $\pi$-Hall subgroups of $S$, i.e., $G^\varphi$ is a
subgroup of $\Sym_{k_\pi(S)}$ acting without stable points.

If $2\not\in\pi$, then by \cite[Theorem~A]{GroConjOddOrder}, $S\in
C_\pi$, so $G^\varphi\simeq\Sym_1$ has a stable point and this case is impossible.

If $3\not\in\pi$ then \cite[Corollary~5.3]{RevVdoDpiFinal} implies that $S$ always possesses an $\Aut(S)$-in\-va\-ri\-ant class of conjugate $\pi$-Hall
subgroups, so again $G$ has a stable point on $K_\pi(S)$ and this case is impossible.

Suppose that $2,3\in\pi$. Then
\cite[Theorem~1.1]{RevVdoArxive} implies that either $k_\pi(S)\in\{1,2,3,4\}$, or $k_\pi(S)=9$. As in the case $2\not\in\pi$ we see, that
$k_\pi(S)=1$ is impossible. By \cite[Lemma~8.2]{RevVdoArxive} we obtain that if $k_\pi(S)=9$, then $S$ possesses an $\Aut(S)$-invariant
conjugacy class of $\pi$-Hall subgroup, hence $G$ has a stable point on $K_\pi(S)$. So we may assume that $2\le k_\pi(S)\le 4$.

Assume that $k_\pi(S)=2$ and $G\not\in E_\pi$. Then $G^\varphi\simeq \Sym_2$, so there exists an element $x\in G$ such  that
$x^\varphi=(1,2)$. Clearly we may assume that the order of $x$ is a power of $2$. By Lemma \ref{HallExist} we obtain $\langle
x,S\rangle\not\in E_\pi$.

Assume that $k_\pi(S)=3$. Since $G$ has no fixed points on $K_\pi(G)$, it acts transitively on $K_\pi(S)$, so there exists an
element $x\in G$ such  that $x^\varphi=(1,2,3)$.
Clearly we may assume that the order of $x$ is a power of $3$. By Lemma \ref{HallExist} we obtain~${\langle x,S\rangle\not\in E_\pi}$.

Assume finally that  $k_\pi(S)=4$. Then  $G^\varphi$, as a subgroup of $\Sym_4$, acts without fixed points. Now every
subgroup of $\mathrm{Sym}_4$ that does not fix a point
contains an element of type $(i,j)(k,l)$. Consider an element $x\in G$ in the preimage  of such an element  $(i,j)(k,l)$ under the
homomorphism $G\rightarrow
\mathrm{Sym}_4$.
Since $\vert (i,j)(k,l)\vert=2$, we can choose
$x$ of order a power of $2$. Furthermore $x$  stabilizes no class of conjugate $\pi$-Hall subgroups of $S$,
hence by
Proposition \ref{base} and Lemma \ref{HallExist}, the group $\langle x,S\rangle$ has no
$\pi$-Hall subgroups.\qed
\medskip 

\noindent {\bfseries Note.} Using the classification of $\pi$-Hall subgroups it is possible to show that if $k_\pi(S)=3$, then $S$
possesses an $\Aut(S)$-in\-va\-ri\-ant class of conjugate $\pi$-Hall subgroups. So the element $x$ from Lemma \ref{EpiCyclic} may be chosen to be a
$2$-ele\-ment. However this fact needs a rather complicated analysis of simple $E_\pi$-groups with $3$ conjugacy classes of
{$\pi$-Hall} subgroup and we do not give it here.
\medskip

\noindent {\bfseries Proof of Theorem \ref{ExistenceCriterion}.} Assume by contradiction that  $G$ is
an $E_\pi$-gro\-up, with composition series $1=G_0<G_1<G_2<\ldots<G_k=G$, and that $\Aut_G(G_i/G_{i-1})\not\in
E_\pi$ for some $i$. If
$G_i/G_{i-1}$ is Abelian,
then $\Aut(G_i/G_{i-1})$ is cyclic. Since $$\Aut_G(G_i/G_{i-1})\leq
\Aut(G_i/G_{i-1}),$$ we obtain that  $\Aut_G(G_i/G_{i-1})\in E_\pi$. Thus $G_i/G_{i-1}$ is non-abelian. By
Lemma \ref{EpiCyclic} we have that
$2,3\in\pi\cap\pi(G)$ and there exists $x\in \Aut_G(G_i/G_{i-1})$ such that $\langle x,G_i/G_{i-1}\rangle\not\in E_\pi$ and the order of
$x$ is either a power of $2$, or a power of $3$. Let $y\in
N_G(G_i/G_{i-1})$ be a preimage of $x$; it may be chosen to be a
$2$-ele\-ment or a $3$-ele\-ment. By Sylow's theorem there exists a $\pi$-Hall subgroup $H$ of
$G$ such that $y\in H$. By Proposition \ref{base}, $(H\cap G_i)G_{i-1}/G_{i-1}$ is
a $\pi$-Hall subgroup of $G_i/G_{i-1}$. Moreover in view of the choice of $y$
we have  $y\in N_G((H\cap G_i)G_{i-1}/G_{i-1})$, hence $x$ (as the image of $y$ in $ N_G(G_i/G_{i-1})/C_G(G_i/G_{i-1})$) normalizes
a $\pi$-Hall subgroup $(H\cap G_i)G_{i-1}/G_{i-1}$ of $G_i/G_{i-1}$. By Lemma
\ref{HallExist} we obtain that $\langle x,G_i/G_{i-1}\rangle\in E_\pi$, a
contradiction.\qed
\medskip

Now we are able to prove the corollaries.\medskip

\noindent{\bfseries Proof of Corollary \ref{SubDirectProduct}.}
It is enough to prove that if $G$ is a finite group with normal subgroups $M,N$
such that $M\cap N=1$ and $G/M,G/N\in E_\pi$, then $G\in E_\pi$. Clearly we may assume that both $M$ and $N$ are proper subgroups of $G$.
Let $1=G_0<G_1<\ldots<G_n=G$  be a composition series of $G$ which is a refinement of
a chief series of $G$ through $M$. If $G_i>M$, then the condition
$G/M\in E_\pi$ and Theorem \ref{ExistenceCriterion} imply that
$\Aut_G(G_i/G_{i-1})\in E_\pi$. If $G_i\leq M$, then the condition $M\cap N=1$
and Lemma \ref{InducedAutomorphism} imply that $\Aut_G(G_i/G_{i-1})\simeq \Aut_{G/N}((G_iN/N)/(G_{i-1}N/N))$. Condition $G/N\in E_\pi$ and
Theorem \ref{ExistenceCriterion} imply that $ \Aut_{G/N}((G_iN/N)/(G_{i-1}N/N))\in E_\pi$. Therefore by
Corollary \ref{GrossExistenceCriterionImproved} we obtain that~${G\in
E_\pi}$.\qed\medskip

\noindent{\bfseries Proof of Corollary \ref{HomomorphicImage}.} Let $G\in E_\pi$, $A\unlhd G$ and $M/A$ be a $\pi$-Hall subgroup of $G/A$.
Let $1=G_0<G_1<\ldots<G_k=G$ be a composition series of $G$ which is a
refinement of a chief series of $G$ through $A$, so $G_n=A$ for some $n$. Let
$1=M_n/A<M_{n+1}/A<\ldots<M_m/A=M/A$ be a composition series which is a refinement of
a chief series of $M/A$. Then
$$1=G_0<G_1<\ldots<G_n=A=M_n< M_{n+1}<\ldots
< M_m=M$$ is a composition series of $M$ which is a refinement of a chief
series of $M$. Since  $M/A$ is a $\pi$-gro\-up, then,  for each $i>n$, $\Aut_M(M_i/M_{i-1})$ satisfies $E_\pi$. For every
non-Abelian composition factor $G_i/G_{i-1}$ with $i\le n$ we have that
$$G_i/G_{i-1}\leq \Aut_M(G_i/G_{i-1})\leq
\Aut_G(G_i/G_{i-1}).$$ By Theorem  \ref{ExistenceCriterion},
$\Aut_G(G_i/G_{i-1})$ satisfies $E_\pi$ and 
\cite[Corollary~3.3]{GrossExistence} implies that $\Aut_M(G_i/G_{i-1})$ satisfies $E_\pi$. By Theorem \ref{GrossExistenceCriterion}, $M$
satisfies~$E_\pi$.  Hence there exists a $\pi$-Hall subgroup $H$ of $M$ and $H$ is a $\pi$-Hall subgroup of $G$. Since $M/A$ is a
$\pi$-gro\-up, then~${M/A=HA/A}$. \qed



\end{document}